\documentclass{lmcs}
\pdfoutput=1

\usepackage[utf8]{inputenc}

\usepackage{lastpage}
\lmcsdoi{17}{1}{2}
\lmcsheading{}{\pageref{LastPage}}{}{}%
{Apr.~02,~2018}{Jan.~22,~2021}{}

\usepackage{url,amsmath,amsthm,amssymb,xspace,hyperref}

\mathchardef\hyphen="2D
\newcommand{\Th}{\mathrm{Th}}

\newcommand{\LFP}{\mathrm{LFP}}
\newcommand{\FO}{\mathrm{FO}}

\newcommand{\FOP}[1]{\FO\hyphen#1}
\newcommand{\LFPP}[1]{\LFP\hyphen#1}

\newcommand{\OP}{\mathrm{OP}}
\newcommand{\NOP}{\mathrm{NOP}}
\newcommand{\IP}{\mathrm{IP}}
\newcommand{\NIP}{\mathrm{NIP}}
\newcommand{\SOP}{\mathrm{SOP}}
\newcommand{\NSOP}{\mathrm{NSOP}}
\newcommand{\TP}{\mathrm{TP}}

\newcommand{\TPO}{\mathrm{TP}_1}

\newcommand{\TPT}{\mathrm{TP}_2}
\newcommand{\NTPT}{\mathrm{NTP}_2}

\newcommand{\RG}{\mathrm{RG}}
\newcommand{\PEQ}{\mathrm{PEQ}}
\newcommand{\GIS}{\mathrm{GIS}}

\newcommand{\TINF}{T_\infty}
\newcommand{\TRG}{T_\mathrm{rg}}
\newcommand{\DLO}{\mathrm{DLO}}
\newcommand{\TFEQ}{T_{\mathrm{feq}}}
\newcommand{\TORG}{T_{\mathrm{org}}}
\newcommand{\TABA}{T_{\mathrm{aba}}}

\newcommand{\A}{\mathbf{A}}
\newcommand{\M}{\mathbf{M}}

\newcommand{\bit}{\mathrm{bit}}
\newcommand{\factor}{\mathrm{factor}}
\newcommand{\fin}{\mathrm{fin}}

\newcommand{\Fraisse}{Fra\"iss\'e\xspace}

\theoremstyle{thmC}
\newtheorem{factC}[thm]{Fact}

\begin{document}

\title{Tameness in least fixed-point logic and McColm's conjecture}

\begin{abstract}
We investigate four model-theoretic tameness properties in the context of least fixed-point logic over a family of finite structures. We find that each of these properties depends only on the elementary (i.e., first-order) limit theory, and we completely determine the valid entailments among them. In contrast to the context of first-order logic on arbitrary structures, the order property and independence property are equivalent in this setting.

McColm conjectured that least fixed-point definability collapses to first-order definability exactly when proficiency fails. McColm's conjecture is known to be false in general. However, we show that McColm's conjecture is true for any family of finite structures whose limit theory is model-theoretically tame.
\end{abstract}

\author{Siddharth Bhaskar}
\address{Department of Computer Science, University of Copenhagen, Denmark}
\email{sbhaskar@di.ku.dk}

\author{Alex Kruckman}
\address{Department of Mathematics and Computer Science, Wesleyan University, USA}
\email{akruckman@wesleyan.edu}


\maketitle

\section{Introduction}

Least fixed-point ($\LFP$) logic is obtained by extending first-order ($\FO$) logic by a quantifier denoting the least fixed-point of a relational operator. The difference between $\FO$ and $\LFP$ definability over classes of finite structures is a central question in finite model theory. McColm conjectured that the existence of arbitrarily long elementary inductions would suffice to separate $\LFP$ from $\FO$ \cite{McC}. This conjecture was refuted by two separate constructions due to Gurevich, Immerman, and Shelah \cite{GIS}. However, instances of this conjecture remain interesting: in particular, any resolution of the \emph{ordered conjecture}, which states that $\LFP$ is more expressive than $\FO$ over every class of totally ordered structures, would resolve a longstanding open problem in computational complexity \cite{KV}.

Recent work, e.g., Adler and Adler \cite{AA}, has shown that certain model-theoretic tameness properties introduced by Shelah \cite{S} are relevant to finite model theory, generalizing assumptions like \emph{bounded cliquewidth} and \emph{bounded treewidth} which permit, e.g., fast algorithms for formula evaluation \cite{Cour}. Here we study four such properties, $\NOP$, $\NIP$, $\NSOP$, and $\NTPT$. We show that any counterexample to McColm's conjecture must fail all of them, thus establishing in a precise sense that any such counterexample must be complicated (Theorem \ref{tame-McColm}).

This result was foreshadowed by McColm himself, who observed in~\cite{McC} that any counterexample to his conjecture must fail another model-theoretic tameness property, viz., the negation of the finite cover property ($\mathrm{NFCP}$). Our work complements Lindell and Weinstein \cite{LW00}, who show that any counterexample cannot be \emph{recursion-theoretically} tame.

In the course of this investigation, we formulate $\LFP$ versions of these four properties. In $\FO$ logic, $\NOP\Rightarrow\NSOP$ and $\NOP\Rightarrow\NIP\Rightarrow\NTPT$, but no other entailments hold in general between these four properties. By contrast, in $\LFP$ logic over finite structures, we have $\NOP\Leftrightarrow\NIP\Rightarrow\NTPT\Rightarrow\NSOP$ (Corollary \ref{OP=IP}). Moreover, both of the implications in $\NIP\Rightarrow\NTPT\Rightarrow\NSOP$ are strict (Theorem \ref{thm:otherway}). We find the equivalence $\NOP\Leftrightarrow\NIP$ especially remarkable: intuitively, it says that ``order implies randomness'' in this context.

Finally, we find that each of the properties of families of finite structures that we study, viz., $\FO = \LFP$, proficiency, and $\LFP\hyphen(\NOP, \NIP, \NSOP, \NTPT)$, depend only on the elementary limit theory of the class of structures (Lemma~\ref{lem:Steve's-2}, Corollary~\ref{cor:Steve's-1-cor}, and Corollary \ref{FO theory defns of LFP props}). This relatively innocuous observation seems not to have been explicitly mentioned before, but the resulting shift in perspective, from classifying structures to classifying theories, brings these questions closer to the spirit of classical model theory.

\subsection{Least fixed-point logic}

We assume familiarity with $\FO$ and $\LFP$ definability, and we very briefly review the latter. (See, e.g., Libkin \cite{Lib} for a reference.) Ordinary (first-order) variables are denoted by lowercase Latin letters, e.g., $x,y,z$. Relational (second-order) variables are denoted by uppercase Latin letters, e.g., $P,Q,R,S$. Every relational variable comes with an arity, but this is not made explicit in the notation.

Following the typical model-theoretic convention, we will also use, e.g., $x$, to denote a tuple of (first-order) variables, not just a single variable. We write $|x|$ to indicate the length of the tuple $x$. Below, when we write, e.g., $\varphi(x,S)$, we mean that $\varphi$ is a formula, $x$ is a tuple of first-order variables (which includes all free first-order variables in $\varphi$), and $S$ is a single second-order variable which is free in $\varphi$ (but $\varphi$ may have other free second-order variables).

We use boldface letters, especially $\A$ and $\M$, to denote structures. Their respective domains are denoted in lightface, e.g., $A$ and $M$. We denote a class of structures by $\mathcal{C}$; its elements are always assumed to share a common signature.

\begin{defi}
A formula $\varphi(x,S)$ is \emph{positive elementary in $S$} if each occurrence of $S$ is in the scope of an even number of negations. In addition, $\varphi(x,S)$ is \emph{operative} if it is positive elementary in $S$ and the arity of $S$ is $|x|$.
\end{defi}

An operative formula $\varphi(x,S)$ is so called because for any structure $\A$, given together with  interpretations of all free second-order variables in $\varphi$ other than $S$, $\varphi$ defines a monotone operator $A^{|x|} \to A^{|x|}$ by
\[ R \mapsto \{a\in A^{|x|} \mid \A \models \varphi(a,R)\}.\]
For every ordinal $\alpha$, we define
\[ \varphi^\alpha = \begin{cases}
\emptyset & \text{if } \alpha =0 \\
\{a\in A^{|x|} \mid \A \models \varphi(a,\varphi^\beta) \} & \text{if } \alpha = \beta + 1 \\
\bigcup_{\beta < \alpha}\varphi^\beta & \text{if } \alpha \text{ is a limit}
\end{cases}\]
The relations $\varphi^\alpha$ are called the \emph{stages} of $\varphi$ on $\A$. The \emph{closure ordinal} $\|\varphi \|_\A$ is the least ordinal $\Gamma$ such that $\varphi^\Gamma = \varphi^{\Gamma +1}$. (Note that if $A$ is finite, $\|\varphi\|_\A$ must be finite.) The relation $\varphi^\Gamma$ is the least fixed-point of the monotone operator defined by $\varphi$, and is written $\varphi^\infty$.

\begin{defi}\label{def:lfp}
The class of \emph{$\LFP$ formulas} is obtained from the class of atomic formulas (in first- and second-order variables) by closing under boolean connectives, first-order quantifiers, and the least fixed-point quantifier: If $\varphi(x,S)$ is an operative $\LFP$ formula and $t$ is a tuple of terms of length $|x|$, then  $[\mathbf{lfp}\,Sx.\varphi](t)$ is an $\LFP$ formula, in which the free first-order variables are those in $t$ and the free second-order variables are those in $\varphi$, except for $S$, which is bound by the quantifier.

The standard semantics of first-order logic are extended to the least fixed-point quantifier as follows: A structure $\A$ (given together with an interpretation of the free first- and second-order variables) satisfies  $[\mathbf{lfp}\,Sx.\varphi](t)$ if and only if the interpretation of the tuple of terms $t$ is in the relation $\varphi^\infty$.
\end{defi}

\begin{defi}
A \emph{query} $R(x)$ of arity $n$ over a class of structures $\mathcal{C}$ is an isomorphism-invariant family of $n$-ary relations $R^A \subseteq A^n$ for each $\A \in \mathcal{C}$. It is \emph{$\LFP$-definable} in case it is defined by some $\LFP$ formula with no free second-order variables, uniformly over all structures in $\mathcal{C}$. An important special case is a \emph{boolean} query, whose arity is zero. $\LFP$-definable boolean queries are defined by $\LFP$ sentences.

We say $\LFP = \FO$ over a class of structures $\mathcal{C}$ if every query which is $\LFP$-definable over $\mathcal{C}$ is defined by a first-order formula. Otherwise, we say $\LFP>\FO$ over $\mathcal{C}$.
\end{defi}

\begin{defi}
An operative formula $\varphi(x,S)$ is \emph{basic} if it is first-order (i.e., it does not contain any instances of the least fixed-point quantifier) and it has no free second-order variables other than $S$.
\end{defi}

\begin{rem}\label{rem:normalform}
We are primarily concerned with definability over families of finite structures. Immerman \cite{I}, building on the work of Moschovakis \cite{Mos}, proves the following normal form for $\LFP$ formulas over finite structures: \[(Qy)\,([\mathbf{lfp}\,Sx.\varphi](t))\] for a basic operative formula $\varphi$, a tuple of terms $t$, and a string of first-order quantifiers $Qy$ binding some of the free variables in $t$.

When working with a class of finite structures, this normal form allows us to restrict attention to $\LFP$ formulas containing a single least fixed-point quantifier. In particular, we only need to consider least fixed-point quantification of basic operative formulas. On the other hand, when we want to show that a particular relation is $\LFP$-definable, we will freely make use of the full syntax in Definition~\ref{def:lfp}.
\end{rem}

\begin{defi}
Operative formulas $\varphi(x,S)$ and $\psi(y,T)$ are \emph{complementary on finite structures} in case $|x| = |y|$ and for every finite structure $\A$, \[\A \models \forall z\, \left([\mathbf{lfp}\,Sx.\varphi](z)\leftrightarrow \lnot [\mathbf{lfp}\,Ty.\psi](z)\right).\]
\end{defi}

\begin{factC}[\cite{I}]
For every basic operative formula $\varphi(x,S)$, there exists a basic operative formula $\psi(y,T)$, such that $\varphi(x,S)$ and $\psi(y,T)$ are complimentary on finite structures.
\end{factC}

\subsection{Proficiency and McColm's conjecture}

\begin{defi}
We say that a class of finite structures $\mathcal{C}$ is \emph{proficient} if there exists a basic operative formula $\varphi(x,S)$ such that
\[ \sup \{ \| \varphi \|_\A : \A \in \mathcal{C} \} = \omega. \]
\end{defi}

\begin{rem}
\label{LFP=FO} For any basic operative formula $\varphi(x,S)$ and any finite $n$, the stage $\varphi^n$ is definable by a first-order formula, uniformly over all structures. It follows immediately that if a class of finite structures $\mathcal{C}$ is not proficient, then $\LFP = \FO$ over $\mathcal{C}$. This observation was made originally in
\cite{McC}.
\end{rem}

In 1990, McColm~\cite{McC} conjectured that the following
three properties are equivalent, for any family of finite structures $\mathcal{C}$: 
\begin{enumerate}
\item $\mathcal{C}$ is proficient. 
\item $\FO < \LFP$ over $\mathcal{C}$. 
\item $\FO < \mathcal{L}_{\infty\omega}^{\omega}$ over $\mathcal{C}$. 
\end{enumerate}
The implication from 2 to 1 is Remark~\ref{LFP=FO}, and the implication from 3 to 1 is also easy. In 1992, Kolaitis and Vardi \cite{KV} proved the equivalence
of 1 and 3. In 1994, Gurevich, Immerman, and Shelah \cite{GIS}
constructed two examples of a proficient family of structures for
which $\FO = \LFP$, thus establishing that 1 does not imply 2.

Historically, the equivalence of 1 and 2 has been called \emph{McColm's first conjecture} and the equivalence of 1 and 3 \emph{McColm's second conjecture} (see, e.g., \cite{KV}). In the interests of brevity, we will simply say \emph{McColm's conjecture} to mean the equivalence of 1 and 2, following the usage in \cite{GIS}. 

\subsection{Elementary limit theories}

In model theory, the most important invariant of a structure is its theory. In the present paper, the most important invariant of a family of finite structures is its limit theory.

\begin{defi} Let $\mathcal{C}$ be a class of finite structures. For a first-order sentence $\varphi$, we write $\mathcal{C} \models \varphi$ in case all but finitely many structures in $\mathcal{C}$ satisfy $\varphi$. The \emph{(elementary) limit theory} of a family $\mathcal{C}$ of finite structures is \[\Th(\mathcal{C}) = \{\varphi\mid \varphi\text{ is a first-order sentence, and }\mathcal{C}\models \varphi\}.\]
\end{defi}

Unlike the theory $\Th(\A)$ of a structure $\A$, limit theories are not always complete --- nor consistent! But it is easy to see that $\Th(\mathcal{C})$ is consistent if and only if $\mathcal{C}$ is infinite.
Henceforth, we will only consider infinite families of finite structures. 

\begin{lem}\label{lem:closure}
For any class of finite structures $\mathcal{C}$, $\Th(\mathcal{C})$ is closed under logical consequence. 
\end{lem}
\begin{proof}
Suppose $\Th(\mathcal{C})\models \varphi$. By compactness,  $\Delta\models \varphi$ for some finite subset $\Delta\subseteq \Th(\mathcal{C})$. Each sentence in $\Delta$ is true in all but finitely many structures in $\mathcal{C}$, and $\Delta$ is finite, so all but finitely many structures in $\mathcal{C}$ satisfy $\Delta$, and hence satisfy $\varphi$. Thus $\varphi\in \Th(\mathcal{C})$. 
\end{proof}

\begin{defi}
A first-order theory $T$ has the \emph{finite model property} if for every sentence $\varphi$, if $T\models \varphi$, then $\varphi$ has a finite model.
\end{defi}

\begin{lem}\label{lem:fmp}
Let $T$ be a countable first-order theory. Then $T$ has the finite model property if and only if $T \subseteq \Th(\mathcal{C})$ for some infinite class of finite structures $\mathcal{C}$.
\end{lem}
\begin{proof}
Suppose $T$ has the finite model property. Enumerate $T$ as $\{\varphi_n\mid n\in \mathbb{N}\}$, and let $\psi_n = \bigwedge_{i=0}^n \varphi_i$ for each $n\in \mathbb{N}$. Then $T\models \psi_n$, so $\psi_n$ has a finite model $\A_n$ (which we may assume to be distinct from $\A_m$ for all $m<n$). Letting $\mathcal{C} = \{A_n\mid n\in \mathbb{N}\}$, each sentence $\varphi_n\in T$ is satisfied by $\A_m$ for all $m\geq n$, so $\mathcal{C}\models \varphi_n$, and hence $T\subseteq \Th(\mathcal{C})$.

Conversely, suppose $T\subseteq \Th(\mathcal{C})$ for some infinite class of finite structures $\mathcal{C}$. If $T\models \varphi$, then $\varphi\in \Th(\mathcal{C})$, by Lemma~\ref{lem:closure}. Since $\mathcal{C}$ is infinite, there exists some $\A\in \mathcal{C}$ such that $\A\models \varphi$, so $T$ has the finite model property. 
\end{proof}

\subsection{Model-theoretic dividing lines}\label{sec:dividinglines}
Stability theory originated in the work of Shelah in his program to classify the models of certain complete first-order theories. He discovered a robust division of theories into ``stable'' and ``unstable.'' The former are ``tame'' in the sense that their definable sets are highly structured, in a way that makes it possible to classify their models (under certain additional hypotheses); the latter are ``wild'' in the sense that they interpret combinatorial objects such as infinite linear orders and random graphs, and they necessarily have too many models to admit a nice classification.

A critical observation about the stable/unstable dichotomy is that stability can be defined by the absence of a simple combinatorial configuration in 
definable sets,
namely the \emph{order property} (described below). Subsequent work in ``neo-stability'' has generalized stability theory to increasingly more inclusive notions of tameness, each of which is defined by the absence of a particular configuration in 
definable sets.
The goal is to find robust dividing lines, such that it is possible to prove structure theorems on the tame side and non-structure theorems on the wild side. For an interactive guide to the these dividing lines, see~\cite{C}; see Hodges \cite{Hodges2} for a discussion of structure theorems.

Any formula $\varphi(x;y)$, whose free variables are partitioned into a tuple $x$ and a tuple $y$, defines a bipartite graph relation $R_\varphi$ between $M^{|x|}$ and $M^{|y|}$ for any structure $\M$, by $(a,b)\in R_\varphi$ if and only if $\M\models \varphi(a;b)$. From another point of view, such a formula defines a family of $S_\varphi$ of subsets of $M^{|x|}$: writing $\varphi(M;b)$ for $\{a\in M^{|x|}\mid \M\models \varphi(a;b)\}$, we let $S_\varphi = \{\varphi(M;b)\mid b\in M^{|y|}\}$. By a combinatorial configuration, we usually mean some concrete property of the graph $R_\varphi$ or the family of sets $S_\varphi$.

We will now give the precise definitions of the combinatorial properties we will consider in this paper, all
of which are originally due to Shelah~\cite{S}. See below for further discussion.

\begin{defi}\label{def:properties}
Let $\varphi(x;y)$ be any formula ($\FO$ or $\LFP$), whose free variables are partitioned into a tuple $x$ and a tuple $y$. Let $n$ be a natural number, and let $\M$ be a structure. 
\end{defi}
\begin{itemize}
\item $\varphi$ has an \emph{$n$-instance of the order property} ($\OP(n)$) in $\M$ if there exist tuples $a_1,\dots,a_n\in M^{|x|}$ and $b_1,\dots,b_n\in M^{|y|}$ such that $\M\models \varphi(a_i;b_j)$ if and only if $i\leq j$. 

\item $\varphi$ has an \emph{$n$-instance of the independence property} ($\IP(n)$) in $\M$ if there exist tuples $a_i\in M^{|x|}$ for all $i\in \{1,\dots,n\}$ and $b_X\in M^{|y|}$ for all $X\subseteq \{1,\dots,n\}$ such that $\M\models \varphi(a_i;b_X)$ if and only if $i\in X$. 

\item $\varphi$ has an \emph{$n$-instance of the strict order property} ($\SOP(n)$) in $\M$ if there exist tuples $b_1,\dots,b_n\in M^{|y|}$ such that $\varphi(M;b_i)\subseteq \varphi(M;b_j)$ if and only if $i\leq j$.

\item We say $\varphi(x;y)$ has an \emph{$n$-instance of the tree property
of the second kind} ($\TPT(n)$) in $\M$ if there are tuples $b_{i,j}\in M^{|y|}$
for $1\le i,j\le n$ such that for any $i$ and any $j\neq k$,  $\varphi(M;b_{i,j})\cap \varphi(M;b_{i,k}) = \emptyset$, but for any function $f\colon\{1,\dots,n\}\to\{1,\dots,n\}$,
\[
\bigcap_{i=1}^n \varphi(M;b_{i,f(i)}) \neq \emptyset.
\]
\end{itemize}

\begin{defi}
Let $\mathcal{C}$ be a class of structures. For each property $P$ in $\{\OP,\IP,\SOP,\TPT\}$, we say $\mathcal{C}$ \emph{has} $\LFPP{P}$ (resp.\ $\FOP{P}$) if there exists an $\LFP$ formula (resp.\ an $\FO$ formula) $\varphi(x;y)$ such that for each $n$, there exists a structure $\M\in \mathcal{C}$ such that $\varphi$ has $P(n)$ in $\M$. 

Let $T$ be a theory. We say that $T$ has $\FOP{P}$ if the class of models of $T$ has $\FOP{P}$.

If a class of structures $\mathcal{C}$ or a theory $T$ does not have $(\FO/\LFP)\hyphen\OP$ (resp.\ $\IP$, $\SOP$, $\TPT$), we say it has $(\FO/\LFP)\hyphen\NOP$ (resp.\ $\NIP$, $\NSOP$, $\NTPT$).
\end{defi}

\begin{rem}
Our definitions of these properties differ from the standard definitions in model theory, which consist of a single infinite configuration, rather than a sequence of finite configurations. For example, according to the standard definition, a theory $T$ has the order property if there exists an $\FO$ formula $\varphi(x;y)$, a model $\M\models T$, and tuples $(a_n)_{n\in \mathbb{N}}$ in $M^{|x|}$ and $(b_n)_{n\in \mathbb{N}}$ in $M^{|y|}$ such that $\M\models \varphi(a_i,b_j)$ if and only if $i\leq j$. 

In the context of a first-order theory $T$, our definitions are equivalent to the standard ones, by an application of the compactness theorem. But compactness is not available in the context of $\LFP$ definability, and the standard infinitary definitions are not meaningful over classes of finite structures. 

Another advantage of using finite configurations is that the presence or absence of an $n$-instance of one of our properties in a structure $\M$ is expressible by a single sentence. For example, $\varphi(x;y)$ has $\OP(n)$ in $\M$ if and only if \[\M \models \exists x_1\dots \exists x_n \exists y_1\dots \exists y_n \left(\bigwedge_{i\leq j} \varphi(x_i;y_j)\land \bigwedge_{i>j} \lnot \varphi(x_i;y_j)\right).\] For $P\in \{\OP,\IP,\SOP,\TPT\}$, we denote by $P_\varphi(n)$ the sentence expressing that $\varphi(x;y)$ has $P(n)$. This leads immediately to the following lemma.
\end{rem}

\begin{lem}\label{lem:FO-P passage}
Let $\mathcal{C}$ be a class of finite structures. For any $P$ in $\{\OP, \IP, \SOP, \TPT\}$, $\mathcal{C}$ has $\FOP{P}$ if and only if $\Th(\mathcal{C})$ has $\FOP{P}$. 
\end{lem}
\begin{proof}
Suppose $\mathcal{C}$ does not have $\FOP{P}$. Then for every $\FO$ formula $\varphi(x;y)$, there is some $n\in \mathbb{N}$ such that for all $\M\in \mathcal{C}$, $\varphi(x;y)$ does not have $P(n)$ in $\M$. That is, $\lnot P_\varphi(n)$ is true in every structure in $\mathcal{C}$, so $\lnot P_\varphi(n)\in \Th(\mathcal{C})$. Thus $\varphi(x;y)$ does not have $P(n)$ in any model of $T$, so $T$ does not have $\FOP{P}$. 

Conversely, suppose $\Th(\mathcal{C})$ does not have $\FOP{P}$. Then for every $\FO$ formula $\varphi(x;y)$, there is some $n\in \mathbb{N}$ such that for all $\M\models \Th(\mathcal{C})$, $\varphi(x;y)$ does not have $P(n)$ in $\M$. By Lemma~\ref{lem:closure}, $\lnot P_\varphi(n)\in \Th(\mathcal{C})$, so $\varphi(x;y)$ does not have $P(n)$ in any structures in $\mathcal{C}$, except for finitely many exceptions. Each of these exceptional structures are finite, so there is some maximum $N$ such that $\varphi(x;y)$ has $P(N)$ in any structure in $\mathcal{C}$. Thus $\mathcal{C}$ does not have $\FOP{P}$.
\end{proof}

Intuitively, a formula $\varphi(x;y)$ has the order property if arbitrarily long linear orders are represented in the bipartite graph $R_\varphi$, in the sense that the ``half-graphs'' appear as induced subgraphs: $a_1,\dots,a_n$ and $b_1,\dots,b_n$ with $a_i R_\varphi b_j$ if and only if $i\leq j$. The independence property and the strict order property are two natural strengthenings of this condition: $\IP$ is equivalent to the condition that \emph{arbitrary} bipartite graphs appear as induced subgraphs of $R_\varphi$, and $\SOP$ says that arbitrarily long linear orders appear as chains in the family of sets $S_\varphi$. 

It is not hard to see that any formula with $\IP(n)$ or $\SOP(n+1)$ in a structure has $\OP(n)$ in that structure. At the level of of complete first-order theories, the converse is true (but the same formula need not serve as the witness). This important dichotomy is due to Shelah: an unstable ($\OP$) theory must exhibit order ($\SOP$) or randomness ($\IP$). 

\begin{factC}[{\cite[Theorem II.4.7]{S}}]\label{fact:dichotomy}
A first-order theory $T$ has $\FOP{\OP}$ if and only if it has $\FOP{\IP}$ or $\FOP{\SOP}$.
\end{factC}

The tree property of the second kind is admittedly somewhat less intuitive. Roughly speaking, a formula $\varphi(x;y)$ has the tree property of the second kind if the family of sets $S_\varphi$ includes arbitrarily many arbitrarily large families of disjoint sets, such that these families interact independently. It is not hard to see that any formula with $\TPT(n)$ in a structure has $\IP(n)$ in that structure.

The name $\TPT$ comes from another important dichotomy identified by Shelah: a theory is called \emph{simple} if it does not have the tree property ($\TP$), and a theory has the tree property if and only if it has the tree property of the first kind ($\TPO$) or the tree property of the second kind ($\TPT$). Unlike $\TPT$, the configurations defining the properties $\TP$ and $\TPO$ are visibly related to trees. We will not consider the properties $\TP$ and $\TPO$ in this paper.

\subsection{Examples via \Fraisse limits}

\Fraisse theory is a fruitful source of examples in model theory and provides an important connection between classes of finite structures and the model theory of complete first-order theories. If a class $\mathcal{C}$ of finite structures is isomorphism-closed and countable up to isomorphism, and has the hereditary property, the joint embedding property, and the amalgamation property, then it admits a unique countable \emph{\Fraisse limit} $\M_\mathcal{C}$. This structure is universal and homogeneous for $\mathcal{C}$, in the sense that a finite structure is in $\mathcal{C}$ if and only if it embeds in $\M_\mathcal{C}$, and any two such embeddings are conjugate by an automorphism of $\M_\mathcal{C}$. The theory $\Th(\M_\mathcal{C})$ is called the \emph{generic theory} of $\mathcal{C}$. For more on \Fraisse theory, see~\cite[Section 7.1]{Hodges}.

Here are some examples of generic theories, and which properties from Definition~\ref{def:properties} they do and do not satisfy:

\begin{itemize}
    \item $\TINF$, the theory of an infinite set with no additional structure. This is the generic theory of the class of finite sets. It has $\FOP{\NOP}$ (and hence $\FOP{\NIP}$, $\FOP{\NSOP}$, and $\FOP{\NTPT}$).
    
    \item $\DLO$, the theory of dense linear orders without endpoints. This is the generic theory of the class of finite linear orders. It has $\FOP{\SOP}$ (and hence $\FOP{\OP}$), but $\FOP{\NIP}$ (and hence $\FOP{\NTPT}$).
    
    \item $\TRG$, the theory of the random graph. This is (by definition) the generic theory of the class of finite graphs. It has $\FOP{\IP}$ (and hence $\FOP{\OP}$), but $\FOP{\NSOP}$ and $\FOP{\NTPT}$.

   \item $\TORG$, the generic theory of the class of finite graphs equipped with an ordering of their vertices. It has $\FOP{\IP}$ and $\FOP{\SOP}$ (and hence $\FOP{\OP}$), but $\FOP{\NTPT}$.
   
   \item $\TFEQ$, the generic theory of the class of finite parameterized equivalence relations. The language consists of two unary predicates $O$ and $P$ (for ``objects'' and ``parameters''), and one ternary relation symbol $E$. A parameterized equivalence relation is a structure such that $O$ and $P$ partition the domain, and for every element $a$ satisfying $P$, the binary relation $E(a,x,y)$ is an equivalence relation on the elements satisfying $O$. This theory has $\FOP{\TPT}$ (and hence $\FOP{\IP}$ and $\FOP{\OP}$) but $\FOP{\NSOP}$.
   
   \item $\TABA$, the theory of atomless Boolean algebras. This is the generic theory of the class of finite Boolean algebras. It has $\FOP{\SOP}$ and $\FOP{\TPT}$ (and hence $\FOP{\IP}$ and $\FOP{\OP}$).
\end{itemize}

We will return to several of these examples in Section~\ref{sec:entailments} below.

\section{\texorpdfstring{$\FO$}{FO} characterizations of \texorpdfstring{$\LFP$}{LFP} properties}\label{sec:Proficiency-sOP}
Many important properties of a class of finite structures $\mathcal{C}$ depend only on the elementary limit theory $\Th(\mathcal{C})$, in the sense that any two classes of structures with the same limit theory agree on the property in question. Lemma \ref{lem:FO-P passage} shows this holds of the FO properties from Definition \ref{def:properties}. In this section, we show that it additionally holds for proficiency, ``FO = LFP over $\mathcal{C}$,'' and all the LFP properties from Definition \ref{def:properties}.

\subsection{Proficiency and \texorpdfstring{$\FO = \LFP$}{FO = LFP}}

By Remark~\ref{LFP=FO}, each finite stage $\varphi^n$ of a basic operative formula $\varphi(x,S)$ is definable by a first order formula, uniformly over all structures. We also use $\varphi^n(x)$ to denote any such formula. The context that the symbol $\varphi^n$ appears in will distinguish whether we mean the formula or the relation it defines.

\begin{defi}
A theory $T$ is \emph{not proficient} in case for every basic operative formula $\varphi(x,S)$, there is a natural number $n$ such that \[ T \models (\forall x)(\varphi^n(x) \leftrightarrow \varphi^{n+1}(x)).\]
Otherwise, we say that $T$ is \emph{proficient}.
\end{defi}

\begin{lem}\label{lem:Steve's-2}
Let $\mathcal{C}$ be a class of finite structures. Then $\mathcal{C}$ is proficient if and only if $\mathrm{Th}(\mathcal{C})$ is proficient.
\end{lem}

\begin{proof}
Suppose that $\mathcal{C}$ is not proficient. Then for each basic operative formula $\varphi(x,S)$, there is a bound $n\in \mathbb{N}$ such that $\|\varphi\|_\A \leq n$ for all $\A\in \mathcal{C}$. Therefore, for all $\A \in \mathcal{C}$, $ \A \models (\forall x)(\varphi^n(x)\leftrightarrow \varphi^{n+1}(x)),$ and hence $\mathrm{Th}(\mathcal{C}) \models (\forall x)(\varphi^n(x)\leftrightarrow \varphi^{n+1}(x)).$

Conversely, suppose that $\mathrm{Th}(\mathcal{C})$ is not proficient, and let $\varphi$ be any basic operative formula. Then, there is some $n \in \mathbb{N}$ such that $ \mathrm{Th}(\mathcal{C}) \models (\forall x)(\varphi^n(x)\leftrightarrow \varphi^{n+1}(x)).$ Hence, for all but finitely many $\A \in \mathcal{C}$, $\A \models (\forall x)(\varphi^n(x)\leftrightarrow \varphi^{n+1}(x)),$ and therefore $\| \varphi \|_\A \le n$ for all but finitely many $\A \in \mathcal{C}$. Since there are only finitely many exceptional structures, $\sup \{ \|\varphi \|_\A : \A \in \mathcal{C} \}$ must still be finite.
\end{proof}

\begin{cor}
Let $\mathcal{C}$ and $\mathcal{D}$ be classes of finite structures. 
If $\Th(\mathcal{C}) = \Th(\mathcal{D})$, then $\mathcal{C}$ is proficient if and only if $\mathcal{D}$ is proficient.
\end{cor}

The proof of Lemma~\ref{lem:Steve's-2} actually shows that if $\Th(\mathcal{C})\subseteq \Th(\mathcal{D})$, and $\mathcal{D}$ is proficient, then $\mathcal{C}$ is proficient. Similar refinements, replacing equality of limit theories with containment, can be observed for many of the results in this section.

\begin{defi}
Let $\varphi(x,S)$ be a basic operative formula, and let $\mathcal{C}$ be a class of finite structures. We say that \emph{$\varphi^\infty$ is elementary over $\mathcal{C}$} if there is a first-order formula $\gamma(x)$ which defines the query $\varphi^\infty$ over $\mathcal{C}$. 
\end{defi}

\begin{lem}\label{lem:fixingfinite}
$\varphi^\infty$ is elementary over $\mathcal{C}$ if and only if there is a first-order formula $\gamma(x)$ which defines the query $\varphi^\infty$ over all but finitely many structures in $\mathcal{C}$. 
\end{lem}
\begin{proof}
One direction is trivial. For the other direction, suppose $\gamma(x)$ defines $\varphi^\infty$ over all but finitely many structures in $\mathcal{C}$. Since every finite structure is determined up to isomorphism by a first-order sentence, and every automorphism-invariant relation on a finite structure is definable by a first-order formula, we can modify $\gamma(x)$ so that it defines $\varphi^\infty$ in each of the finitely many exceptional cases. 
\end{proof}

Note that the sentence \[\forall x\,( [\mathbf{lfp}\,Sx.\varphi](x) \leftrightarrow \gamma(x)),\] which expresses that $\gamma$ defines $\varphi^\infty$, is not first-order. So it does not follow directly from Lemma~\ref{lem:fixingfinite} that elementarity of $\varphi^\infty$ over $\mathcal{C}$ is a property of the limit theory $\Th(\mathcal{C})$. Nevertheless, this turns out to be true, as we will now show. 

\begin{lem}\label{lem:Steve's-1}
Let $\varphi(x,S)$ be a basic operative formula. Then $\varphi^\infty$ is elementary over $\mathcal{C}$ if and only if there exists a first-order formula $\theta(x)$ such that
\begin{align*}
\forall x\, (\varphi(x,\theta)\leftrightarrow \theta)&\in \Th(\mathcal{C}), \text{and}\\
 \forall x\,( \psi(x,\neg \theta) \leftrightarrow \neg \theta) &\in \Th(\mathcal{C}),
\end{align*}
where $\psi(x,S)$ is a basic operative formula which is complementary to $\varphi(x,S)$ on finite structures.
\end{lem}

\begin{proof}
Suppose $\varphi^\infty$ is elementary relative to $\mathcal{C}$, witnessed by $\theta$. In an arbitrary structure $\A\in \mathcal{C}$, $\theta$ defines $\varphi^\infty$, which is a fixed-point for $\varphi$. Since $\psi(x,S)$ and $\varphi(x,S)$ are complementary on finite structures, $\neg \theta$ defines $\psi^\infty$, which is a fixed-point for $\psi$. So $\A$ satisfies the sentences in the statement of the lemma. 

Conversely, suppose these two sentences are in $\Th(\mathcal{C})$. Then for all but finitely many structures $\A\in \mathcal{C}$, the relation $\theta^\A$ defined by $\theta$ over $\A$ is a fixed-point of $\varphi$, and its complement $(\lnot\theta)^\A$, which is defined by $\lnot \theta$, is a fixed point of $\psi$. Since $\varphi^\infty$ and $\psi^\infty$ are the least fixed-points of $\varphi$ and $\psi$, $\varphi^\infty \subseteq \theta^\A$ and $\psi^\infty \subseteq (\neg \theta)^\A$. Since $\varphi^\infty$ and $\psi^\infty$ are complements, $\varphi^\infty = \theta^\A$. This is true for all but finitely many structures in $\mathcal{C}$, so $\varphi^\infty$ is elementary over $\mathcal{C}$ by Lemma~\ref{lem:fixingfinite}.
\end{proof}

\begin{cor}\label{cor:Steve's-1-cor}
If $\Th(\mathcal{C}) = \Th(\mathcal{D})$, then $\LFP=\FO$ over $\mathcal{C}$ if and only if $\LFP=\FO$ over~$\mathcal{D}$.
\end{cor}

\begin{proof}
Suppose that every $\LFP$-definable query over $\mathcal{C}$ is $\FO$-definable. Consider an arbitrary $\LFP$-definable query $R$ over $\mathcal{D}$. By the normal form for $\LFP$ formulas over finite structures (Remark~\ref{rem:normalform}), we may assume that $R$ is defined by the $\LFP$ formula $(Qy)\,([\mathbf{lfp}\,Sx.\varphi](t))$. Then it suffices to show that $\varphi^\infty$ is elementary over $\mathcal{D}$.

By assumption, $\varphi^\infty$ is elementary over $\mathcal{C}$. But $\Th(\mathcal{C}) = \Th(\mathcal{D})$, so by Lemma~\ref{lem:Steve's-1}, $\varphi^\infty$ is elementary over $\mathcal{D}$. The converse follows in the same way. 
\end{proof}

\begin{rem}
Even though $\FO = \LFP$ is a property of limit theories, we have \emph{not} proven (and in fact it is not true) that  $\FO = \LFP$ over $\mathcal{C}$ if and only if $\FO = \LFP$ over the class of models of $\Th(\mathcal{C})$. This stands in contrast to the properties FO-$P$ and proficiency, which do pass (Lemmata \ref{lem:FO-P passage} and \ref{lem:Steve's-2}) between $\mathcal{C}$ and models of $\mathrm{Th}(\mathcal{C})$.
\end{rem}

We learned Lemmata  \ref{lem:Steve's-2} and \ref{lem:Steve's-1} from Steven Lindell through personal communication. As an immediate consequence of Lemma~\ref{lem:Steve's-2} and Corollary~\ref{cor:Steve's-1-cor},  we deduce the following.

\begin{cor}
If $\Th(\mathcal{C}) = \Th(\mathcal{D})$, then  $\mathcal{C}$ satisfies McColm's conjecture if and only if $\mathcal{D}$ does.
\end{cor}

\begin{rem}
If we were to define the \emph{$\LFP$ limit theory} to be the set of all $\LFP$ sentences that hold of all but finitely many structures in $\mathcal{C}$, then both $\FO = \LFP$ and proficiency could easily seen to depend only on the $\LFP$ limit theory. One might naturally wonder whether the $\LFP$ limit theory itself depends only on the (elementary) limit theory; this would contain Lemma~\ref{lem:Steve's-2} and Corollary~\ref{cor:Steve's-1-cor} as special cases.

However, this is not the case. For example, the family of even-sized linear orders and the family of odd-sized linear orders are two families with the same (elementary) limit theory: the complete theory of infinite discrete linear orders with endpoints. But they are distinguished by their $\LFP$ limit theories, since parity of the domain is an $\LFP$-definable boolean query over ordered structures. 
\end{rem}

\subsection{Model-theoretic dividing lines}

We will now show that each $\LFP$ property defined in Definition \ref{def:properties} depends only on the elementary limit theory of a family of finite structures. We start by identifying  proficiency with $\LFPP{\SOP}$. Key to this argument is the $\LFP$-definability of the \emph{stage comparison relation} \cite{Mos}.

\begin{defi}
For any basic operative formula $\varphi(x,S)$ with $|x| = k$, any structure $\A$, and any $a\in A^k$, the \emph{stage of $a$},  $\|a\|_\varphi$, is the least ordinal $\alpha$ such that $a\notin \varphi^\alpha$, or $\infty$ if $\alpha\notin \varphi^\infty$. The \emph{stage comparison relation} $\preceq_\varphi$ is defined by $a \preceq_\varphi b$ if and only if $\|a\|_\varphi \leq \|b\|_\varphi$. 
\end{defi}

\begin{factC}[{\cite[Theorem 2A.2]{Mos}}]\label{fact:stagecomp}
For any basic operative formula $\varphi(x,S)$, the stage comparison relation $\preceq_\varphi$ is $\LFP$-definable over the class of all structures.
\end{factC}

By a \emph{partial preorder}, we mean a reflexive transitive relation $\preceq$ on some set $X$. We get a partial order if we take the quotient by the equivalence relation $x \preceq y \wedge y \preceq x$. A partial preorder is \emph{linear} if the associated partial order is linear. By a chain in a partial preorder $\preceq$, we mean a subset $X$ which is linearly ordered by $\preceq$. In particular, for $x,y\in X$, $x\preceq y$ and $y\preceq x$ implies $x = y$. 

Note that for any basic operative formula $\varphi(x,S)$ with $|x| = k$ and any structure $\A$, the stage comparison relation $\preceq_\varphi$ is always a linear preorder on $A^k$, whose associated linear order is a well-order.

\begin{thm}\label{thm:preorder}
\label{main} Let $\mathcal{C}$ be a class of finite structures. The following are equivalent:
\begin{enumerate}
\item $\mathcal{C}$ is proficient.

\item There is some $n \in \mathbb{N}$ and some $\LFP$ formula $\psi$ such that $\psi$ defines a linear preorder on $n$-tuples in every structure in $\mathcal{C}$, and this linear preorder has arbitrarily long finite chains in structures in $\mathcal{C}$. 

\item There is some $n \in \mathbb{N}$ and some $\LFP$ formula $\psi$ such that $\psi$ defines a partial preorder on $n$-tuples in every structure in $\mathcal{C}$, and this partial preorder has arbitrarily long finite chains in structures in $\mathcal{C}$.
\end{enumerate}
\end{thm}
\begin{proof}

($1\Rightarrow 2$) Suppose $\varphi(x,S)$ is a basic operative formula with $|x| = n$, which witnesses proficiency of $\mathcal{C}$. Its stage comparison relation $\preceq_\varphi$, which is $\LFP$-definable by Fact~\ref{fact:stagecomp}, linearly preorders the $n$-tuples from each $\M\in\mathcal{C}$, and this linear preorder contains a chain of length $\|\varphi\|_\M$. By proficiency, there is no finite bound on the lengths of these chains. 

($2\Rightarrow 3$) Trivially.

($3\Rightarrow1$) Suppose $\lambda(y_{1};y_{2})$ defines a partial preorder $\preceq$ which has arbitrarily long finite chains in structures in $\mathcal{C}$. In any finite partial
preorder, we define the \emph{height} of an element to be one more
than the maximum height among its (strict) predecessors, or 0 if it has none.
Since $\lambda$ has arbitrarily long chains, elements in the preorder
defined by $\lambda$ will have arbitrarily large heights.

Let $\varphi(y;T)$ say that all of $y$'s strict predecessors are
in $T$. In symbols, 
\[
\varphi(y;T)\equiv\forall y'\,((\lambda(y',y)\wedge\neg\lambda(y,y'))\to T(y')).
\]
(Notice $T$ occurs positively in $\varphi$.) Then it is easy to
show by induction that the stages $\varphi^n$ of $\varphi$
are exactly those elements of height $<n$, and hence $\varphi$ witnesses that $\mathcal{C}$ is
proficient. 
\end{proof}

\begin{thm}\label{LFP-sOP=proficient}
$\mathcal{C}$ is proficient if and only if it has $\LFPP{\SOP}$.
\end{thm}
\begin{proof}

Suppose that the $\LFP$ formula $\varphi$ witnesses that  $\mathcal{C}$ has $\LFPP{\SOP}$. Define
\[
\psi(y_{1};y_{2})\equiv\forall x\ (\varphi(x;y_{1})\to\varphi(x;y_{2})),
\]
so that for each $\M\in\mathcal{C}$ and $b_1,b_2 \in M^{|y|}$, 
\[
\M\models\psi(b_1;b_2)\text{ if and only if }\varphi(M;b_1)\subseteq\varphi(M;b_2),
\]
where (as in Section~\ref{sec:dividinglines}),  $\varphi(M;b)=\{a \in M^{|x|} \mid \M \models \varphi(a;b)\}$. Then $\psi$ defines a partial preorder on $M^{|y|}$. If $\varphi$ has $\SOP(n)$ in $M$, then this partial preorder on $M^{|y|}$ contains a chain of length $n$. Since
$\varphi$ has the strict order property, we have arbitrarily long chains in structures in $\mathcal{C}$, so $\mathcal{C}$ is proficient by Theorem~\ref{thm:preorder}.

Conversely, suppose that $\mathcal{C}$ is proficient. By Theorem~\ref{thm:preorder}, there exists an $\LFP$ formula $\lambda(y_{1};y_{2})$ that defines a partial preorder with
arbitrarily long chains in structures in $\mathcal{C}$. The formula $\lambda$ itself
witnesses the strict order property: given $n$, pick $\M\in\mathcal{C}$
which contains a chain $b_{1},\dots,b_{n}$. Then 
\[
\lambda(M^{|y_{1}|};b_{1})\subsetneq\dots\subsetneq\lambda(M^{|y_{1}|};b_{n}).\qedhere
\]
\end{proof}

Among all the $\LFP$ properties from Definition~\ref{def:properties}, the strict order property turns out to be the strongest, in that in entails all the others. This contrasts with the first-order case where, in general, the strict order property does not imply the independence property or the tree property of the second kind.

\begin{lem}\label{fundamental-implications}
\label{N_functions}If $\mathcal{C}$ has $\LFPP{\SOP}$, then it also
has $\LFPP{\OP}$, $\LFPP{\IP}$, and $\LFPP{\TPT}$.
\end{lem}
\begin{proof}
As noted in Section~\ref{sec:dividinglines}, $\LFPP{\SOP}$ easily implies $\LFPP{\OP}$.
For the other two properties, we consider the family $\mathcal{N}$ of all finite linear orders. Identify the unique linear order of size $n$ with
the set $n=\{0,1,\dots,n-1\}$ equipped with its natural ordering. 
It is well known that over $\mathcal{N}$, the graphs of addition and multiplication are $\LFP$-definable; hence, so is the graph of exponentiation \cite{L}. Therefore, since the relations $\mathrm{bit}(x;y)$:
\[
``\text{the \ensuremath{x}-th bit of \ensuremath{y} base $2$ is $1$}"
\]
and $\mathrm{factor}(x;y,z)$:
\[
``y^{z}\text{ is the largest power of \ensuremath{y} dividing }x"
\]
are first-order definable over $\mathcal{N}$ with addition, multiplication, and exponentiation, they are $\LFP$-definable over $\mathcal{N}$.

The relation $\bit(x;y)$ has $\IP(n)$ in $m$ for sufficiently large $m$, witnessed by $a_{i}=i-1$ for $i\in \{1,\dots,n\}$ and $b_{X}=\sum_{j\in X} 2^{j-1}$ for
$X\subseteq \{1,\dots,n\}$. The relation $\factor(x;y,z)$ has TP2($n$)
in $(m,<)$ for sufficiently large $m$, witnessed by $b_{i,j}=(p_{i},j)$, where $(p_{i})_{i\in\omega}$
is an enumeration of the primes: for any function $f\colon \{1,\dots,n\}\to \{1,\dots,n\}$, we have  $a_{f}=\prod_{i=1}^{n}p_{i}^{f(i)}\in \factor(M;p_i,f(i))$ for all $1\leq i\leq n$. Hence, $\mathcal{N}$ has $\LFPP{\IP}$ and $\LFPP{\TPT}$.

Now suppose $\mathcal{C}$ has $\LFPP{\SOP}$. Since there is some $\LFP$ formula $\psi$ which defines a linear preorder on $n$-tuples with arbitrarily long chains in structures in $\mathcal{C}$ (by Theorem~\ref{main} and Theorem~\ref{LFP-sOP=proficient}), we can repeat the constructions of $\mathrm{bit}$ and $\mathrm{factor}$ above to get formulas witnessing that $\mathcal{C}$ has $\LFPP{\IP}$ and $\LFPP{\TPT}$. 

To be a little more concrete, suppose
$\varphi(x;y)$ witnesses $\LFPP{\IP}$ or $\LFPP{\TPT}$ over $\mathcal{N}$.
Simply replace each variable $v$ in $\varphi$ by $n$ variables $v_{1},\dots,v_{n}$, replace $v = w$ by $(v\leq w\land w\leq v)$, replace $v\leq w$ by $\psi(v_{1},\dots,v_{n};w_{1},\dots,w_{n})$, 
and proceed by induction on the construction of $\varphi$ in the
obvious way. This gives us a new formula $\varphi^{\star}(x^{\star};y^{\star})$
witnessing $\LFPP{\IP}$ or $\LFPP{\TPT}$ over $\mathcal{C}$.
\end{proof}

\begin{thm}\label{FO defns of LFP props}
For any class $\mathcal{C}$ of finite structures and any $P$ in $\{\OP,\IP,\SOP,\TPT\}$, 
$\mathcal{C}$ has $\LFPP{P}$ if and only if $\mathcal{C}$ is proficient or $\mathcal{C}$ has  $\FOP{P}$.
\end{thm}

\begin{proof}
In the forwards direction, if $\mathcal{C}$ has $\LFPP{P}$, but is not proficient, then $\FO = \LFP$ over $\mathcal{C}$ (Remark~\ref{LFP=FO}), so $\mathcal{C}$ has $\FOP{P}$.

Conversely, if $\mathcal{C}$ is proficient, then it has $\LFPP{P}$ by Theorems \ref{LFP-sOP=proficient} and \ref{fundamental-implications}. Otherwise, if $\mathcal{C}$ has $\FOP{P}$, then it trivially has $\LFPP{P}$ as well.
\end{proof}

Theorem \ref{FO defns of LFP props}, when combined with Lemmata~\ref{lem:FO-P passage} and~\ref{lem:Steve's-2}, has the immediate consequence that all the $\LFP$ properties depend only on the elementary limit theory of a family of structures.

\begin{cor}\label{FO theory defns of LFP props}
Suppose $\mathcal{C}$ and $\mathcal{D}$ are families of structures with the same limit theory, and let $P$ be any property in $\{\OP,\IP,\SOP,\TPT\}$. Then
$\mathcal{C}$ has $\LFPP{P}$ if and only if $ \mathcal{D}$ has $\LFPP{P}$.
\end{cor}

It also has the nice consequence that any tame class of structures (in the sense of first-order model theory) satisfies McColm's conjecture.

\begin{thm}\label{tame-McColm}
For any family of finite structures $\mathcal{C}$, if $\mathcal{C}$ has $\FOP{\NOP}$, $\FOP{\NIP}$, $\FOP{\NTPT}$, or $\FOP{\NSOP}$, then $\mathcal{C}$ satisfies McColm's conjecture.
\end{thm}

\begin{proof}
Let $P\in \{\OP,\IP,\SOP,\TPT\}$ be a property such that $\mathcal{C}$ does not have $\FOP{P}$. To show that $\mathcal{C}$ satisfies McColm's conjecture, it suffices to show that if $\mathcal{C}$ is proficient, then $\LFP\neq \FO$ over $\mathcal{C}$. 

So we assume $\mathcal{C}$ is proficient. Then $\mathcal{C}$ has $\LFPP{P}$, by Theorem~\ref{LFP-sOP=proficient} and Lemma~\ref{fundamental-implications}. Since it does not have $\FOP{P}$, $\LFP\neq \FO$ over $\mathcal{C}$, as desired.
\end{proof}

\section{Entailments between \texorpdfstring{$\LFP$}{LFP} properties}\label{sec:entailments}

We continue towards determining all valid entailments among the $\LFP$ properties. First, we show that an important dichotomy remains true in the $\LFP$ context.

\begin{thm}\label{OP=sOP or IP}
For any family $\mathcal{C}$ of finite structures, 
$\mathcal{C}$ has $\LFPP{\OP}$ if and only if $\mathcal{C}$ has $\LFPP{\SOP}$ or $\LFPP{\IP}$.
\end{thm}

\begin{proof}
Since both $\LFPP{\SOP}$ and $\LFPP{\IP}$ entail $\LFPP{\OP}$, it suffices to show the forwards direction.

Suppose $\mathcal{C}$ has $\LFPP{\OP}$. If $\mathcal{C}$ is proficient, then it has $\LFPP{\SOP}$ by Theorem~\ref{LFP-sOP=proficient}. If $\mathcal{C}$ is not proficient, then $\mathcal{C}$ has $\FOP{\NSOP}$  and $\FOP{\OP}$ by Theorem~\ref{FO theory defns of LFP props}, so $\Th(\mathcal{C})$ has $\FOP{\NSOP}$ and $\FOP{\OP}$ by Lemma~\ref{lem:FO-P passage}. It follows that $\Th(\mathcal{C})$ has $\FOP{\IP}$ by Fact~\ref{fact:dichotomy}. Hence $\mathcal{C}$ has $\FOP{\IP}$, and therefore $\LFPP{\IP}$.  
\end{proof}

By Theorems \ref{fundamental-implications}, \ref{OP=sOP or IP}, and propositional reasoning, we obtain the following.

\begin{cor}\label{OP=IP}
For any family $\mathcal{C}$ of finite structures, 
\[ \mathcal{C} \models \LFPP{\OP}\iff 
\mathcal{C} \models \LFPP{\IP}.\]
Therefore,
\[ \mathcal{C} \models \LFPP{\SOP} \implies
\mathcal{C} \models \LFPP{\TPT} \implies
\mathcal{C} \models \LFPP{\OP} \iff 
\mathcal{C} \models \LFPP{\IP} .\]
\end{cor}

We would like to give examples showing that the first two implications above are strict. To do this, we will employ countably categorical first-order theories with the finite model property. 

A theory is \emph{countably categorical} if it 
has only one countable model up to isomorphism.
Various equivalent formulations of countable categoricity were proven in the 50's and 60's by Ryll-Nardzewski, Svenonius, and Engler. These established countable categoricity as a robust and important property of first-order theories. \Fraisse theory is an important source of examples of countably categorical theories: every \Fraisse limit in a finite relational language has a countably category complete theory. For more information and history on countable categoricity, see Chapter 7 of Hodges \cite{Hodges}.

Though seemingly separate notions, proficiency and countable categoricity are intimately related: roughly speaking, non-proficiency is the ``finite variable logic version'' of countable categoricity. The Ryll-Nardzewski theorem~\cite[Theorem 7.3.1]{Hodges} asserts that countable categoricity is equivalent to  the finiteness of the set of complete $n$-types over $T$, for all $n$. Non-proficiency essentially weakens this condition to the finiteness of the set of all $n$-types \emph{in $m$-variable logic}, for all $n$ and $m \ge n$ \cite[Theorem 23]{DLW}. We now prove more carefully that countable categoricity implies non-proficiency, for first-order theories.

\begin{lem}\label{lem:countable categoricity}
Every countably categorical first-order theory is non-proficient.
\end{lem}

\begin{proof}
Suppose that $T$ is countably categorical, fix a basic operative formula $\varphi(x,S)$, and consider the first-order formulas $\varphi^n(x)$ defining its stages. By the Ryll-Nardzewski theorem, there are only finitely many pairwise non-$T$-equivalent formulas with free variables from $x$. Thus there must be some $m \in \mathbb{N}$ and $n < m$ such that $T \models (\forall x)(\varphi^n(x) \leftrightarrow \varphi^m(x))$. Since, for all $j$, $T \models (\forall x)(\varphi^j(x) \rightarrow \varphi^{j+1}(x))$, it must be the case that $T \models  (\forall x)(\varphi^m(x) \leftrightarrow \varphi^{m+1}(x))$, and $\|\varphi\|_\M \leq m$ for all models $\M\models T$.

Since $\varphi$ was chosen arbitrarily, $T$ is non-proficient.
\end{proof}

\begin{rem}\label{rem:ccfmp}
In particular, if $T$ is a countable and countably categorical theory with the finite model property, then there is some family of finite structures $\mathcal{C}$ with limit theory $T$, by Lemma \ref{lem:fmp}. Since $T$ is countably categorical, it is non-proficient, and hence so is $\mathcal{C}$, by Lemma \ref{lem:Steve's-2}. Therefore, FO = LFP over $\mathcal{C}$. In addition, $\mathcal{C}$ inherits any property FO-$P$ (or its negation) from $T$ itself, by Lemma \ref{lem:FO-P passage}.
\end{rem}

To complete the classification, we show that both of the one-way implications in Corollary \ref{OP=IP} are strict.

\begin{thm}\label{thm:otherway}
\label{thm:noncollapse} There exists a class of finite structures
with $\LFPP{\IP}$ but without $\LFPP{\TPT}$ and a class of finite structures
with $\LFPP{\TPT}$ but without $\LFPP{\SOP}$. 
\end{thm}
\begin{proof}

To exhibit a class of structures with a certain combination of LFP-properties, it suffices to exhibit a class of structures with the same combination of FO-properties, over which FO = LFP. By Lemma \ref{lem:countable categoricity} and Remark~\ref{rem:ccfmp}, it suffices to exhibit a countable, complete, and countably categorical theory with the finite model property, with the same combination of FO-properties.
This is exactly what we do. See Section~\ref{sec:dividinglines} for definitions of our example theories. In this proof, we drop the prefix $\FOP{}$. 

For $\IP$ but $\NTPT$, consider $T_{\RG}$, the theory of the random
graph. This is well-known to be countably categorical with the finite
model property, to have $\IP$, and to be simple; a first-order theory
is simple if it does not have the tree property ($\TP$), which implies
that it does not have $\TPT$.

For $\TPT$ but $\NSOP$, consider $\TFEQ$, the generic theory
of parameterized equivalence relations. For discussions of this theory,
see~\cite{CR} and~\cite{K}. In~\cite{CR}, Chernikov and Ramsey
establish (Corollary 6.20) that $\TFEQ$ does not have the property
$\SOP_1$, which implies that it does not have $\SOP$, and (Corollary 6.18)
that $\TFEQ$ is not simple by witnessing $\TPT$ directly.
A proof that $\TFEQ$ has the finite model property is
given in~\cite{K}. 
\end{proof}

We conclude with a list of some simple examples satisfying the various
combinations of properties we have discussed in this paper (see Figure~\ref{fig:examples}).
Even though the $\LFP$ properties in each box in the table are not explicit, we can easily deduce them: in the column $\LFP=\FO$, the $\LFP$ properties agree with the $\FO$ properties, and in the column $\LFP\neq\FO$, each family of structures is proficient, and hence (by Theorem \ref{FO defns of LFP props}) satisfies each of $\LFPP{\SOP}$,  $\LFPP{\TPT}$, $\LFPP{\IP}$, and $\LFPP{\OP}$. 
Since we have established (Lemma~\ref{fundamental-implications}) that there are no classes satisfying
($\NIP$ and $\SOP$ and $\LFP=\FO$) or ($\IP$ and $\NTPT$ and $\SOP$ and $\LFP=\FO$), the table is complete.

\begin{figure}[ht]
\caption{Examples}
\label{fig:examples} \centering

\begin{tabular}{|c||c|c|}
\hline $\FO$ properties & $\LFP=\FO$ & $\LFP\neq \FO$\\
\hline\hline $\NOP$ = ($\NIP$ and $\NSOP$) & $\mathbb{N}_\fin$ & $(\mathbb{N}_\fin,S)$\\
\hline $\NIP$ and $\SOP$ & --- & $(\mathbb{N}_\fin,<)$\\
\hline $\IP$ and $\NTPT$ and $\NSOP$ & $\RG$ & $\RG+(\mathbb{N}_\fin,S)$\\
\hline $\IP$ and $\NTPT$ and $\SOP$ & --- & $\RG+(\mathbb{N}_\fin,<)$\\
\hline $\TPT$ and $\NSOP$ & $\PEQ$ & $\PEQ+(\mathbb{N}_\fin,S)$\\
\hline $\TPT$ and $\SOP$ & $\GIS$ & $\PEQ+(\mathbb{N}_\fin,<)$
\\\hline \end{tabular}
\end{figure}

Here are the definitions of the classes appearing in Figure~\ref{fig:examples}: 
\begin{itemize}
\item $\mathbb{N}_\fin$ is the class of initial segments of $\mathbb{N}$ with no extra
structure. $(\mathbb{N}_\fin,S)$ and $(\mathbb{N}_\fin,<)$ are the classes of structures
with the same domains, but equipped with the successor relation and
the order relation, respectively. 
\item $\RG$ is any class of finite structures with limit theory $T_{\RG}$,
the theory of the random graph. The class of Paley graphs provides an explicit example (see~\cite{BEH}). 
\item $\PEQ$ is any class of finite structures with limit theory
$T_{\mathrm{feq}}^{*}$, the generic theory of parameterized equivalence
relations. Such a class exists by Lemma~\ref{lem:fmp}.  
\item $\GIS$ is any counterexample to McColm's conjecture. For example,
one of the classes of finite structures devised by Gurevich, Immerman,
and Shelah in~\cite{GIS}. 
\item Given classes of finite structures $\mathcal{C}=\{\M_{i}\mid i\in\omega\}$
and $\mathcal{C}'=\{\M'_{i}\mid i\in\omega\}$ in disjoint
languages $L$ and $L'$, respectively, we denote by $\mathcal{C}+\mathcal{C}'$
the family $\{\M_{i}\sqcup \M_{i}'\mid i\in\omega\}$,
where $\M_{i}\sqcup \M_{i}'$ is the disjoint union
of $\M_{i}$ and $\M_{i}'$. We use the fact that
for any property $P$ in $\{\SOP,\TPT,\IP,\OP\}$,
$\mathcal{C}+\mathcal{C}'$ has $\FOP{P}$ if and only if $\mathcal{C}$
has $\FOP{P}$ or $\mathcal{C}'$ has $\FOP{P}$.
\end{itemize}

\section{Further work}
Our results suggest that it may be fruitful to examine $\SOP$, $\TPT$, $\IP$, and $\OP$ beyond the first-order context. In particular, it would be interesting to examine the extend to which weaker fixed-point logics (like transitive closure logic) recover $\TPT$ and $\IP$ from $\SOP$. Another direction is a program to recover complexity-theoretic tameness properties of families of finite structures (like fast formula evaluation) from model-theoretic tameness assumptions, generalizing  assumptions like bounded treewidth and cliquewidth.

In the spirit of classification theory, we might hope to deduce some positive concrete information about, e.g., $\LFPP{\NOP}$ classes of finite structures that distinguish them from the merely stable ($\FOP{\NOP}$) classes. One might hope to develop some kind of asymptotic structure theory (like Shelah's classification of models of certain stable theories) for finite classes which are stable and non-proficient.

Finally, we believe that the observation that $\FO = \LFP$ is a property of the elementary limit theory of a class of finite structures strongly suggests a model-theoretic approach to difficult questions like the ordered conjecture. At the very least, it gives us a new set of powerful tools to test and clarify where, exactly, the difficulty lies.

\section*{Acknowledgements}

We started this project in 2017 when we were both postdocs at Indiana University, Bloomington. We would like to thank Larry Moss and the logic group at IU for their support.

We would also like to acknowledge several people who have read earlier versions of this paper and discussed this material with us. In particular, we thank Cameron Hill for his comments on the first
draft of this paper, Steve Lindell and Scott
Weinstein for being invaluable sources of knowledge of finite model theory, and Greg McColm for inspiring the present line of inquiry. Finally, we are indebted to the anonymous referees, whose comments substantially improved the readability of this paper.

\bibliography{bibliography}

\begin{thebibliography}{DLW96}

\bibitem[AA14]{AA}
Hans Adler and Isolde Adler.
\newblock Interpreting nowhere dense graph classes as a classical notion of
  model theory.
\newblock {\em European J. Combin.}, 36:322--330, 2014.

\bibitem[BEH81]{BEH}
Andreas Blass, Geoffrey Exoo, and Frank Harary.
\newblock Paley graphs satisfy all first-order adjacency axioms.
\newblock {\em J. Graph Theory}, 5(4):435--439, 1981.

\bibitem[Con]{C}
Gabriel Conant.
\newblock Map of the universe.
\newblock \url{https://forkinganddividing.com/}.

\bibitem[Cou90]{Cour}
Bruno Courcelle.
\newblock The monadic second-order logic of graphs. {I}. {R}ecognizable sets of
  finite graphs.
\newblock {\em Inform. and Comput.}, 85(1):12--75, 1990.

\bibitem[CR16]{CR}
Artem Chernikov and Nicholas Ramsey.
\newblock On model-theoretic tree properties.
\newblock {\em J. Math. Log.}, 16(2):1650009, 2016.

\bibitem[DLW96]{DLW}
Anuj Dawar, Steven Lindell, and Scott Weinstein.
\newblock First order logic, fixed point logic and linear order.
\newblock In {\em Computer science logic ({P}aderborn, 1995)}, volume 1092 of
  {\em Lecture Notes in Comput. Sci.}, pages 161--177. Springer, Berlin, 1996.

\bibitem[GIS94]{GIS}
Yuri Gurevich, Neil Immerman, and Saharon Shelah.
\newblock Mccolm's conjecture [positive elementary inductions].
\newblock In Samson Abramsky, editor, {\em Proceedings of the Ninth Annual IEEE
  Symp. on Logic in Computer Science, {LICS} 1994}, pages 10--19. IEEE Computer
  Society Press, July 1994.

\bibitem[Hod87]{Hodges2}
Wilfrid Hodges.
\newblock What is a structure theory?
\newblock {\em Bull. London Math. Soc.}, 19(3):209--237, 1987.

\bibitem[Hod93]{Hodges}
Wilfrid Hodges.
\newblock {\em Model theory}, volume~42 of {\em Encyclopedia of Mathematics and
  its Applications}.
\newblock Cambridge University Press, Cambridge, 1993.

\bibitem[Imm86]{I}
Neil Immerman.
\newblock Relational queries computable in polynomial time.
\newblock {\em Inform. and Control}, 68(1-3):86--104, 1986.

\bibitem[Kru19]{K}
Alex Kruckman.
\newblock Disjoint {$n$}-amalgamation and pseudofinite countably categorical
  theories.
\newblock {\em Notre Dame J. Form. Log.}, 60(1):139--160, 2019.

\bibitem[KV92]{KV}
P.G. Kolaitis and M.Y. Vardi.
\newblock Fixpoint logic vs. infinitary logic in finite-model theory.
\newblock In Andre Scedrov, editor, {\em Proceedings of the Seventh Annual IEEE
  Symp. on Logic in Computer Science, {LICS} 1992}, pages 46--57. IEEE Computer
  Society Press, June 1992.

\bibitem[Lib04]{Lib}
Leonid Libkin.
\newblock {\em Elements of finite model theory}.
\newblock Texts in Theoretical Computer Science. An EATCS Series.
  Springer-Verlag, Berlin, 2004.

\bibitem[Lin]{L}
Steven Lindell.
\newblock Exponentiation is elementarily definable from addition and
  multiplication on finite structures.
\newblock \url{http://ww3.haverford.edu/cmsc/slindell/exponentiation.pdf}.

\bibitem[LW00]{LW00}
Steven Lindell and Scott Weinstein.
\newblock The role of decidability in first order separations over classes of
  finite structures.
\newblock In Martin Abadi, editor, {\em Proceedings of the Fifteenth Annual
  IEEE Symp. on Logic in Computer Science, {LICS} 2000}, pages 45--50. IEEE
  Computer Society Press, June 2000.

\bibitem[McC90]{McC}
Gregory~L. McColm.
\newblock When is arithmetic possible?
\newblock {\em Ann. Pure Appl. Logic}, 50(1):29--51, 1990.

\bibitem[Mos74]{Mos}
Yiannis~N. Moschovakis.
\newblock {\em Elementary induction on abstract structures}.
\newblock North-Holland Publishing Co., Amsterdam-London; American Elsevier
  Publishing Co., Inc., New York, 1974.
\newblock Studies in Logic and the Foundations of Mathematics, Vol. 77.

\bibitem[She90]{S}
S.~Shelah.
\newblock {\em Classification theory and the number of nonisomorphic models},
  volume~92 of {\em Studies in Logic and the Foundations of Mathematics}.
\newblock North-Holland Publishing Co., Amsterdam, second edition, 1990.

\end{thebibliography}
\bibliographystyle{alpha}

\end{document}